\documentclass[sn-mathphys,Numbered]{sn-jnl}

\usepackage{graphicx}%
\usepackage{multirow}%
\usepackage{amsmath,amssymb,amsfonts}%
\usepackage{amsthm}%
\usepackage{mathrsfs}%
\usepackage[title]{appendix}%
\usepackage{xcolor}%
\usepackage{textcomp}%
\usepackage{manyfoot}%
\usepackage{booktabs}%
\usepackage{algorithm}%
\usepackage{algorithmicx}%
\usepackage{algpseudocode}%
\usepackage{listings}%

\usepackage{overpic}
\def\Q{\mathcal{Q}}
\def\C{\mathcal{C}}
\def\S{\mathcal{S}}
\def\R{\mathcal{R}}
\def\Har{\mathcal Har}
\def\Hyp{\mathcal Hyp}
\def\P{\mathbb{P}}

\def\meet{\wedge}

\def\Fa{\colorbox{white}{\small \bf a)}}
\def\Fb{\colorbox{white}{\small \bf b)}}

\theoremstyle{thmstyleone}%
\newtheorem{theorem}{Theorem}
\newtheorem{lemma}{Lemma}%

\begin{document}

\title[Article Title]{Harmonic curves and the beauty of Projective Geometry}

\author[1]{\fnm{José Luis} \sur{Abreu}}\email{
jlabreu@im.unam.mx}

\author*[1]{\fnm{Javier} \sur{Bracho}}\email{jbracho@im.unam.mx}

\affil*[1]{\orgdiv{Instituto de Matemáticas}, \orgname{Universidad Nacional Autónoma de México}, \orgaddress{\street{Ciudad Universitaria}, \city{Mexico City}, \postcode{04510}, 
\country{Mexico}}}

\abstract{The purpose of this paper is to present projective geometry in a synthetic, visual and intuitive style through the central notion of \emph{harmonicity} which leads to \emph{harmonic curves}. This presentation includes new results, unpublished proofs of some classic theorems and a slight reformulation of its axiomatics.}



\maketitle

\section{Introduction} 
A \emph{conic curve} or \emph{section} is defined as the intersection of a \emph{circular cone} with a plane. This concept is inherently metric. However, that conic curves somehow belong to projective geometry is well known, but not clear how. For example, the statement that \emph{ellipses are what we see when we look at a circle}, which is how we should present them to children, makes many of us unconfortable because we didn't hear it in class and don't have a simple proof at hand. Trying to make sense of this incongruent state of affairs, we were lead to define \emph{harmonic curves} in a purely projective style (they are basically the locus of points in the plane that \emph{see} a quadrangle as a harmonic set) and study their properties without any reference to metric concepts. Yes, harmonic curves become conic sections in the euclidean plane, but they have a rich set of purely projective properties so that their study constitutes a beautiful approach to projective geometry. Such approach is what this paper is all about.

Motivated by the technics of perspective drawing, which were discovered by the painters of the renaissance to represent realisticly three dimensional scenes on flat canvases, Girard Desargues (1591-1661) initiated the development of Projective Geometry expanding the concept of space to include \emph{ideal points} (also known as points \emph{at infinity}), \cite{Desargues}. Mostly forgotten through the following two centuries, his visionary work was revived and recognized as fundamental during the first half of the XIX century, when Projective Geometry was firmly established as a field on its own. In 1872 Felix Klein opened his famous \emph{Erlangen Program}, \cite{Erlangen}, with the statement \emph{``Among the advances of the last fifty years in the field of geometry, the development of Projective Geometry occupies the first place''}. 
One of the mathematicians whose work deserved such praise is Karl Georg Christian von Staudt (1798-1867). In his treatise on the subject \cite{VonStaudt}, he proves that the notions of \emph{harmonicity}, i.e. the concepts of harmonic sets and pencils, are entirely projective and independent of metric ones like distances or angles. He also shows that conic curves may be defined through the abstract, and purely projective concept of \emph{polarity}. 

In an attempt to make these facts accesible to high-school teachers, and inspired by \cite{JS} where John Stillwell argues that such goal is not only worthwhile and culturally urgent, but also feasible, the authors of this article developed the \emph{dynamic geometry} system \emph{ProGeo3D} \cite{PG3} which is specialized in projective geometry. 
In particular, it incorporates harmonicity as a construction tool. All the figures of this paper are snapshots of constructions in ProGeo3D. To appreciate their dynamic capacity and the 3D nature of some (e.g., Figure~\ref{Harmonic_theorem}.a), we present them in a web site\footnote{\url{https://arquimedes.matem.unam.mx/harmonic_curves/scenes.html}}. Playing with this system, we discovered several interesting constructions, proofs and new results that do not seem to be included in the known literature on the subject. This paper presents them to propose an alternative approach to projective geometry which is intuitive, synthetic and, in our subjective opinion, beautiful.

We start by reviewing von Staudt's definition of harmonicity, emphasizing the duality it intrinsically carries. That leads us to a simple definition of what we call \emph{harmonic curves} to differentiate them from the classic treatment of conic sections, although in euclidean space they are the same. 
In the following section, we relate them to von Staudt's definition, which uses polarities, but now it must be stated as a theorem, the ``Polarity Theorem'' (Thm.~\ref{teo:Polaridad2D}). Its proof uses the classic idea in projective geometry of going out to 3D, finding a way to work out things there, and coming back to 2D. The first to use our explicit technic out in 3D, was Germinal Pierre Dandelin (1794–1847) to prove Pascal's Hexagon Theorem in \cite{Dandelin}. It is deeply related to ruled surfaces as defined by Hilbert and Cohn-Vossen in \cite{HCV} using only incidence geometry. This suggests the formulation of a new axiom for Projective Geometry, which we call the Equipal Axiom. It is shown to be equivalent to Pappus' Theorem and thus to other equivalent axioms frequently used in projective geometry. In the final section we discuss the axiomatic foundations of our approach in which the use of 3 dimensions is fundamental.

\section{Harmonic sets, pencils and reflections}

One of the seminal contributions of Karl von Staudt was to prove that \emph{harmonicity} (the notion of harmonic conjugates which had been used since antiquity in terms of distances) only depends on incidence using quadrangles.  

A  \emph{quadrangle}, $\Q$, is defined as four points in the projective plane in \emph{general position} (i.e., no three of them are collinear), called its \emph{vertices}, together with 4 lines, called its \emph{sides}, such that their incidence relation is a 4-\emph{cycle} (each object of one type is incident with two of the other); that is, $\Q$ is defined by its vertices (4 points in general position) together with a \emph{dihedral} (i.e., a cyclic but non-oriented) order on them which yields the 4 sides. 

The term ``quadrangle'' is adequate because  at any vertex, its two incident lines (an ``angle'') distinguish  two \emph{adjacent} vertices 
and thus, it also determines the 
\emph{opposite} vertex (as the remaining one); note that \emph{the partition of vertices into opposite pairs determines the quadrangle}. The \emph{center} of the quadrangle $\Q$ is the intersection of its two \emph{diagonal lines} (the lines that join opposite vertices); and its \emph{horizon} is the line joining the intersection of opposite sides, which are called its two \emph{diagonal points}.

The notion of a quadrangle is autodual, but in a  \emph{quadrilateral} the stress is given on the 4 sides, so that its dual-center is the horizon of the corresponding quadrangle whose center is the dual-horizon of the quarilateral. However, we will keep the ``horizon'' and ``center'' terminology: understood always as a line and a point, respectively. These terms are natural because if a square tile is drawn on a canvas, the center of the tile and the horizon of the square tiling of the plane which it generates should be drawn, respectively, at the center and the horizon of the quadrangle which such a drawing determines.

Four points in general position are the vertices of three quadrangles. Their corresponding centers and horizons form its \emph{diagonal triangle}. 

\begin{figure}[H]
         \centerline{
          \href{https://arquimedes.matem.unam.mx/harmonic_curves/viewer.html?view=pg3/Quadrangle_and_harmony.pg3}
     { \begin{overpic}[abs,width=6.1cm]{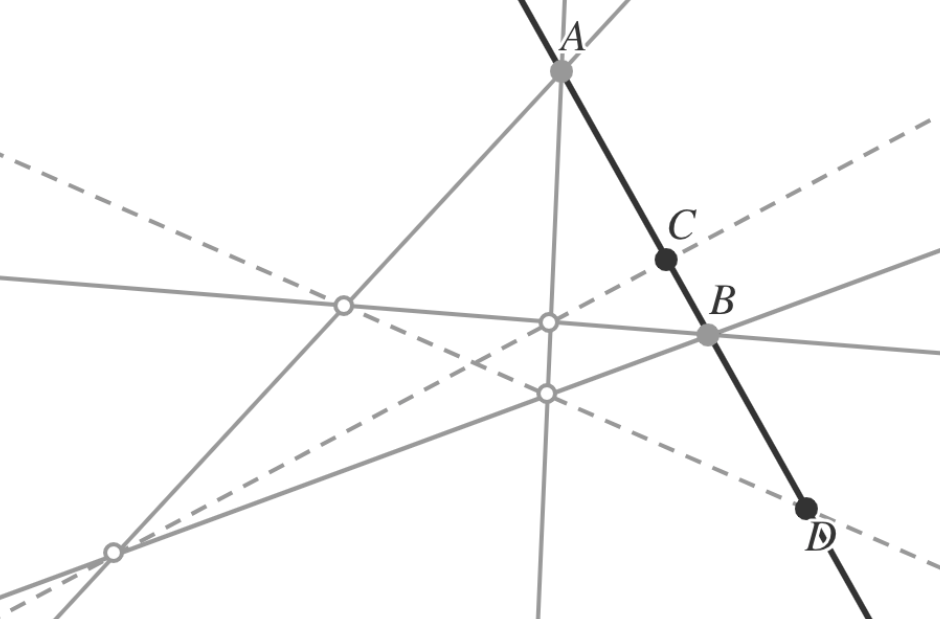}
        \put(-5,0){\Fa}
        \end{overpic} }                 
          \hspace{0.2cm}
         \href{https://arquimedes.matem.unam.mx/harmonic_curves/viewer.html?view=pg3/Harmonic_pencil.pg3} 
         { \begin{overpic}[abs,width=6.1cm]{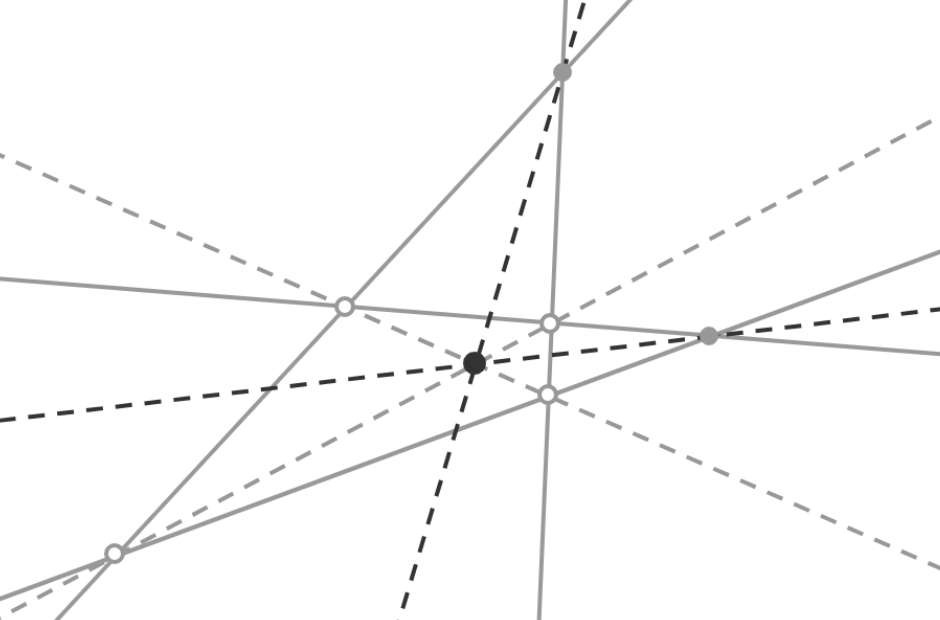}   
           \put(-5,0){\Fb}
         \end{overpic} }}
   \caption{\label{Armonia_y_dual}
    {\bf a)} A harmonic set.  {\bf b)} A harmonic pencil with dashed lines.}  
\end{figure}

Four collinear points  $A, C, B, D$, as in Figure~\ref{Armonia_y_dual}.a, are said to be a \emph{harmonic set},\footnote{The terms ``harmonic quadruple'' or ``harmonic range" are also used, but we stick to ``harmonic set" as in the classic texts \cite{VeblenYoung} and \cite{CoxeterPG}.} if there exists a quadrangle such that its diagonal points are $A$ and $B$ 
(hence the horizon of the quadrangle is their \emph{support line}) 
and the other pair, $C$ and $D$, are incident with the diagonal lines. Dually, four concurrent lines are called a \emph{harmonic pencil}\footnote{The term ``harmonic set of lines'' is also used, e.g. \cite{VeblenYoung, CoxeterPG}; but we will use ``pencil'' for simplicity and to distinguish them immediately from harmonic sets (of points).} if there exists a quadrangle such that one pair of lines are the diagonal lines of the quadrangle 
(hence, its center is the concurrency point of the pencil, 
also called its \emph{center}) and the other two lines are incident to the diagonal points of the quadrangle, see Figure~\ref{Armonia_y_dual}.b. 

As stated, the pairs of elements in the definitions play a different role but, as we will see, they are interchangeable, so that both notions include an explicit dihedral order of the four elements involved, which coincides with their geometric placement (points within a projective line or lines about a point).

Let us now show that these definitions are sound. Given a collinear triple $A, C, B$ with $C$ distinguished, two auxiliary points out of the support line and collinear with $C$ determine a unique quadrilateral $\Q$ as in Figure~\ref{Armonia_y_dual}.a, and therefore produce the point $D$ as the intersection of the other diagonal line with the horizon; this construction, called the \emph{harmonic fourth}, has as outcome the point $D$, called the  \emph{harmonic conjugate} of $C$ with respect to $A$ and $B$. Since for the triple $A, D, B$, one can choose the other opposite pair of vertices of $\Q$ as auxiliary points and then obtain $C$ as outcome, we can further say that  the (unordered) pair of points $C, D$ are  \emph{harmonic conjugates} with respect to $A, B$, \cite{CoxeterPG, VeblenYoung}. 

\begin{theorem}[Harmonic Theorem] The outcome of the harmonic fourth construction does not depend on the choice of the auxiliary points.\end{theorem} 

\begin{figure}[H] \label{Harmonic_theorem}
         \centerline{
          \href{https://arquimedes.matem.unam.mx/harmonic_curves/viewer.html?view=pg3/Harmonic_Theorem_in_3D.pg3}
         { \begin{overpic}[abs,width=6.1cm]
               {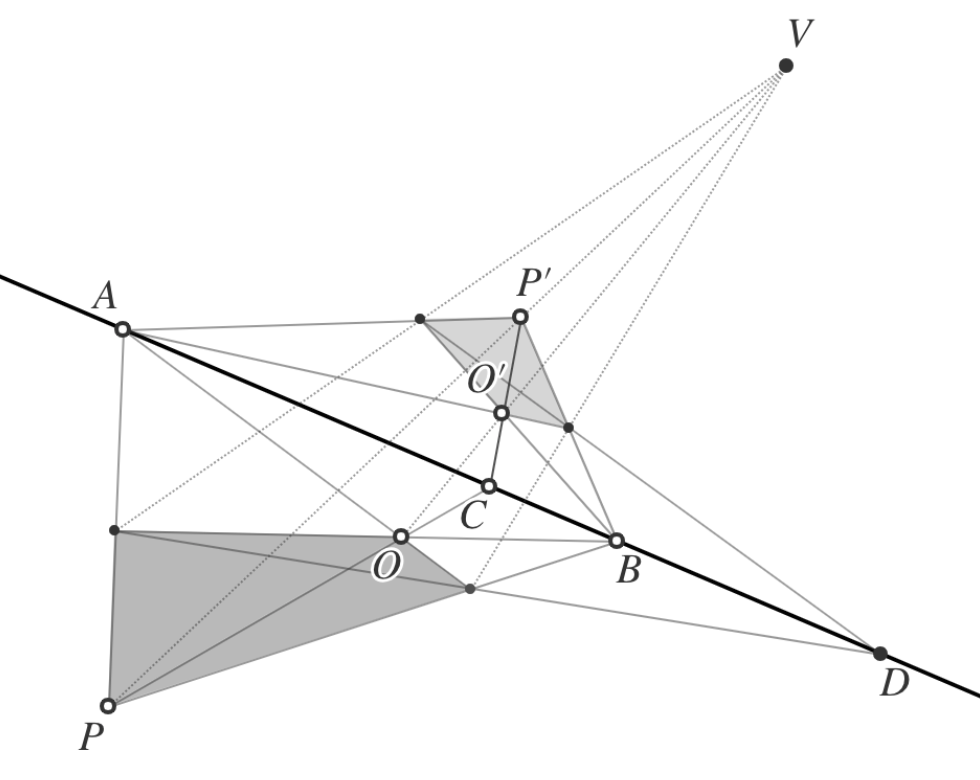}                   
        \put(-5,0){\Fa}
        \end{overpic} }                 
          \hspace{0.2cm}
           \href{https://arquimedes.matem.unam.mx/harmonic_curves/viewer.html?view=pg3/Harmonicity_is_symmetric_with_respect_to_the_pairs_of_points.pg3}{ \begin{overpic}[abs,width=6.1cm]{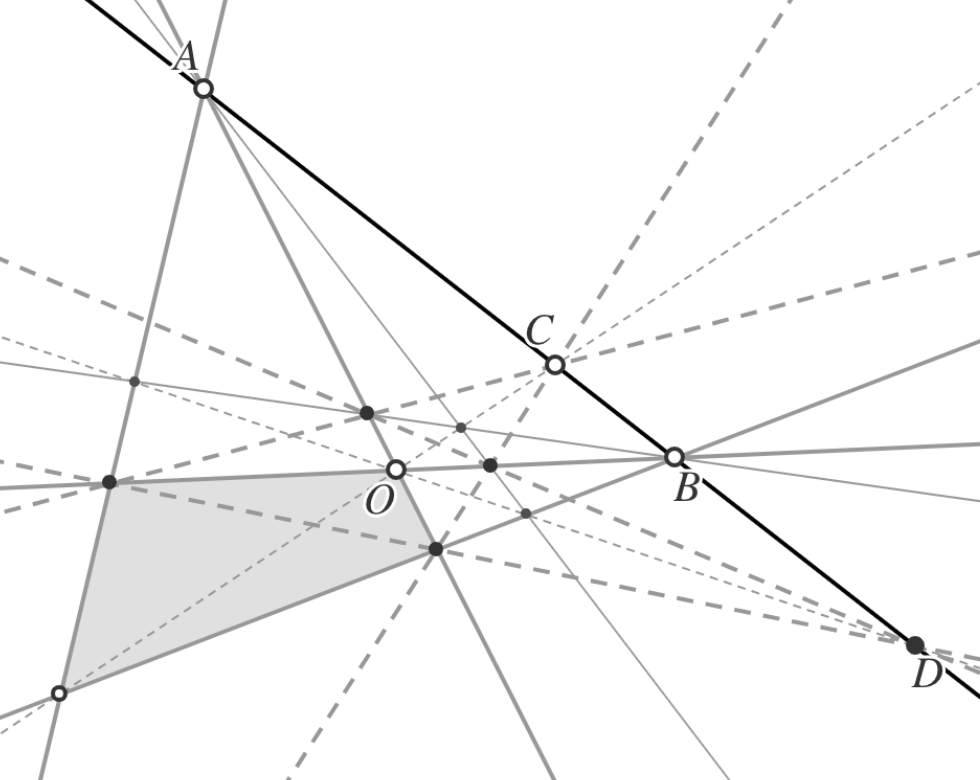}  
           \put(-5,0){\Fb}
         \end{overpic} }}
          \caption{\label{Harmonic_Theorem}
        {\bf a)} Visual proof of the harmonic theorem in 3D. 
                 {\bf b)} Symmetry of harmonic pairs.} 
\end{figure}

The proof is well known and follows from Desargues' Theorem. It may also be directly proved in 3D using only incidence arguments (see Figure~\ref{Harmonic_theorem}.a), and deducing the planar case imediatly from it. We omit the details for brevity.

Finally, to see that the definition of harmonic set is symmetric with regard to the role played by the two pairs of points, extend the quadrangle $\Q$ to a 2 by 2 tiling drawn in perspective, as in Figure~\ref{Harmonic_theorem}.b (the know-how comes from the renaissance artists and the coincidences follow from the Harmonic Theorem). Then, the quadrilateral of (dashed) diagonals not incident with its center $O$ proves that $A, B$ are harmonic conjugates with respect to $C, D$.

Since the point $O$ in Figure~\ref{Harmonic_Theorem}.b may be chosen to be any point not on the support line of the harmonic set, we obtain that any such point \emph{sees} them as harmonic, that is, the lines to them with their dihedral order is a harmonic pencil. Dually, there is also a \emph{harmonic fourth} construction for lines and any line not through the center of a harmonic pencil cuts it in a harmonic set. Thus, \emph{harmonic sets and pencils are preserved by projections}.

The harmonic fourth construction also makes sense in the singular cases when $C=A$ or $C=B$. In them, the outside quadrangle collapses to a line, but the diagonals do not collapse and so the construction holds in the sense of not becoming ambiguous, and it yields $D=A$ or $D=B$, respectively. Therefore, given two (distinct) points $A$ and $B$ in a line $\ell$ we get a well defined map
$$\rho_{A,B}:\ell\to\ell$$
called the \emph{harmonic reflection} of $\ell$ with respect to $A$ and $B$: it fixes these two points and it gives the harmonic conjugate elsewhere. $\rho_{A,B}$ is an involution which interchanges the two segments in which the points $A$ and $B$ divide their projective line. And in particular, it interchanges its ideal point at infinity with the (euclidean) midpoint of $A$ and $B$, making the harmonic fourth construction a very useful tool for perspective drawing.

The natural generalization to the projective plane (space) is the \emph{harmonic reflection} $\rho_{C,m}$ with respect to  a point $C$, called the \emph{center}, and a non-incident line (plane) $m$, called the \emph{mirror}.\footnote{In the plane, Coxeter calls it \emph{harmonic homology} in \cite{CoxeterPG}.}  It is defined on every line $\ell$ through $C$ as the harmonic reflection with respect to $C$ and the intersection of $\ell$ with $m$. If we denote by $\vee$, ``join'', and $\wedge$, ``meet'', the basic projective operations of linear span and intersection, respectivelly, we may write for $X\neq C$:
$$X\cdot\rho_{C,m}=X\cdot\rho_{C,(X\vee C)\wedge m}\,,$$
where note that we write the action of maps or functions on the right.
This notion amalgamates two classic euclidean examples: the central inversions, when the mirror is the line (plane) at infinity, and the reflections when the center is the ideal point in the direction perpendicular to the mirror. 

Harmonic reflections are \emph{collineations} (i.e., they send lines to lines). They act in the dual plane as harmonic reflections in the sense that if $\ell$ is a line different from the mirror $m$, then $\ell, m, {\ell\cdot\rho_{C,m}}, {(\ell\wedge m)\vee C}$ is a harmonic pencil centered at $\ell\wedge m$. 

\begin{lemma}[Klein's Triangle]\label{Lemma_KT}
Given a triangle $ABC$ with respective opposite sides $abc$, then $\{\,id_{\P^2}, \rho_{A,a}, \rho_{B,b}, \rho_{C,c}\,\}$\footnote{$\,\,\P^2$ is the projective plane and $id_{\P^2}$ stands for its identity map.} is the Klein four-group.
\end{lemma}

\begin{figure}[H]
         \centerline{
          \href{https://arquimedes.matem.unam.mx/harmonic_curves/viewer.html?view=pg3/Kleins_triangle_lemma.pg3}{ 
          \includegraphics[width=8.0cm]{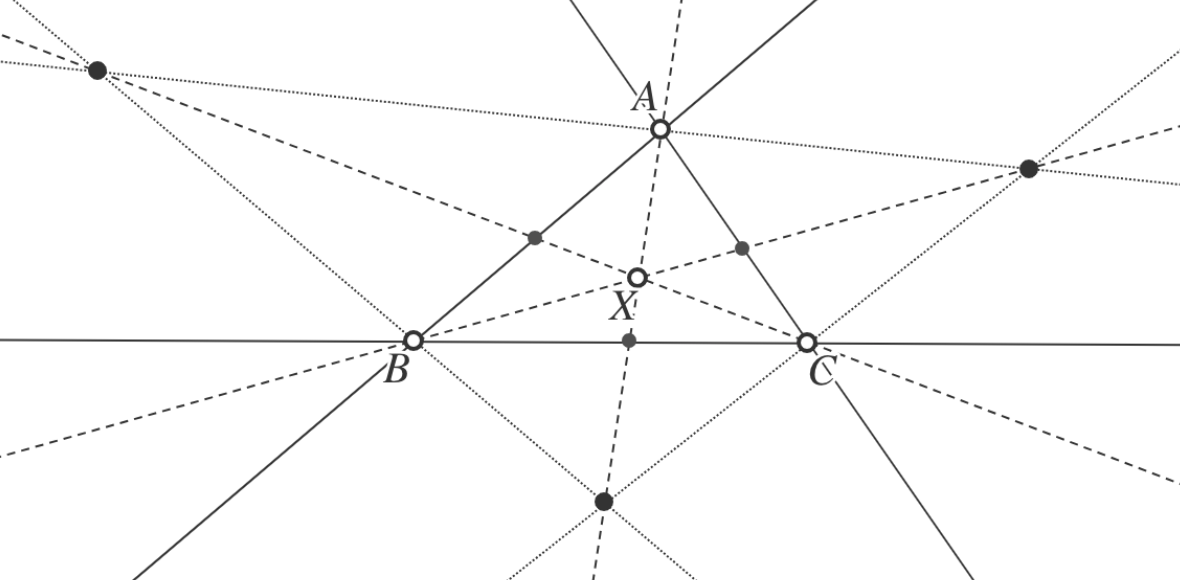}} }
           \caption{\label{Triangulo_de_Klein}
               Klein's Triangle Lemma.}  
\end{figure}

\proof Since the three non-trivial elements are involutions, we need only to show that the composition of any two of them gives the third, which is the definition of the Klein four-group. Consider a point $X$ not in the triangle. We claim that the quadruple $\{\, X, X\cdot\rho_{A,a}, X\cdot\rho_{B,b}, X\cdot\rho_{C,c}\,\}$ has $ABC$ as its diagonal triangle. In Figure~\ref{Triangulo_de_Klein}, the three dashed lines through $X$ have harmonic sets that define the corresponding three points other than $X$. The gray lines from a vertex (say $A$) to one of them (say, $X\cdot\rho_{C,c}$) pass through another one ($X\cdot\rho_{B,b}$) because the two corresponding harmonic sets (in $C\vee X$ and $B\vee X$) are projected to each other from the vertex ($A$) and projections preserve harmonicity. 

It is easy to see that these gray lines through the vertices cut the opposite side in its corresponding harmonic conjugate, and that for points $X$ in the triangle the maps behave as they should. This completes the proof. \qed

Thus, the generic orbits of the Klein four-group associated to a triangle are the quadruples that have it as diagonal triangle, and any of the four triangular regions in which the three lines cut the projective plane are the fundamental regions of the group action which has the vertices as fixed points.

We define the \emph{harmonic group} $\Har(n)$ ($n=1,2,3$) as the group of transformations generated by the harmonic reflections on $\P^n$. Of course, they turn out to be the classic groups of projectivities, but this requires proof.

\section{Harmonic curves and bundles}

We define the \emph{harmonic curve}, $\C_\Q$, of a quadrangle $\Q$ as the locus of points that are centers of harmonic pencils \emph{transversal} to $\Q$. By transversal we mean that each vertex of $\Q$ is incident to a line of the pencil and this correspondence preserves the dihedral orders. Dually, the \emph{harmonic bundle} of a quadrilateral consists of the lines that support a harmonic set transversal to its sides and with corresponding dihedral orders.

\begin{figure}[h]
         \centerline{
         \href{https://arquimedes.matem.unam.mx/harmonic_curves/viewer.html?view=pg3/Quadrangle_with_tangent_lines_to_its_harmonic_curve.pg3}
         { \begin{overpic}[abs,width=6.1cm]{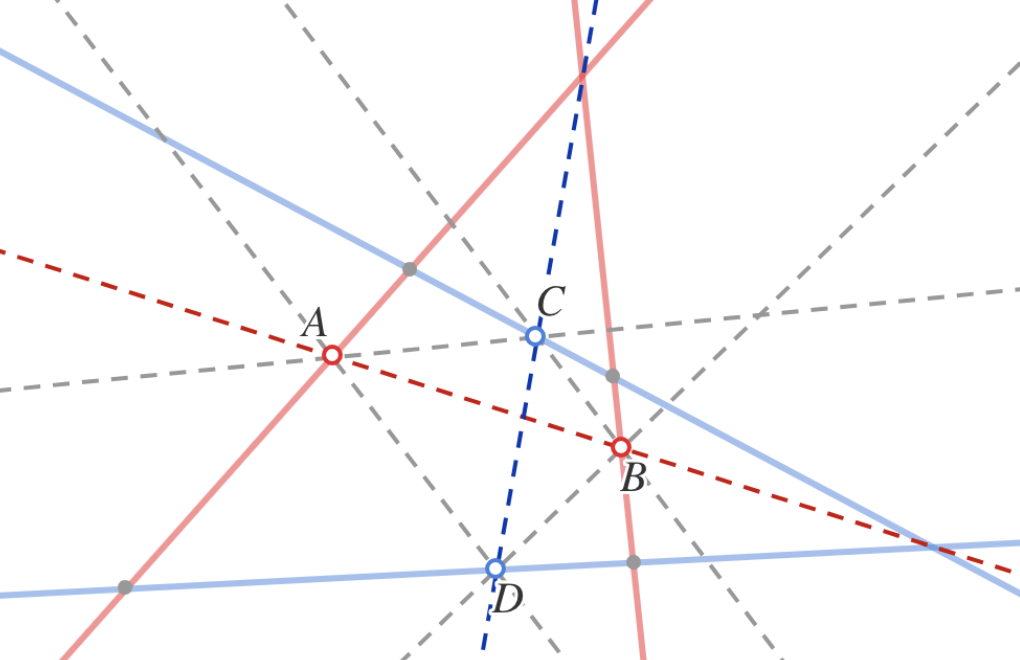}
        \put(-5,0){\Fa}
        \end{overpic} }                 
          \hspace{0.2cm}
          \href{https://arquimedes.matem.unam.mx/harmonic_curves/viewer.html?view=pg3/Harmonic_curve.pg3} 
         { \begin{overpic}[abs,width=6.1cm]{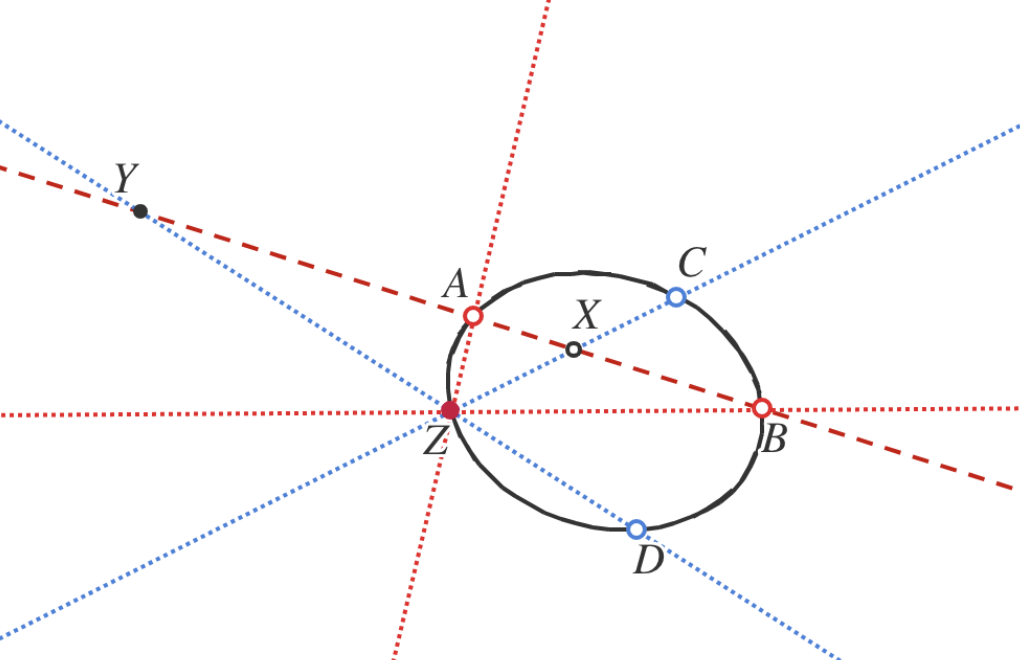}   
           \put(-5,0){\Fb}
         \end{overpic} }}
          \caption{\label{Curva_de_Armonia} {\bf a)} Quadrangle $\Q$ with tangent lines to its harmonic curve $\C_\Q$ at its vertices. {\bf b)} Generic point $Z\in \C_\Q$, where $A,X,B,Y$ is a harmonic set in $A\vee B$.}  
\end{figure}

Consider a quadrangle $\Q$ with vertices $A, C, B, D$. First observe that the vertices are points of its harmonic curve  $\C_\Q$. Indeed, for each vertex, the harmonic conjugate of its diagonal with respect to its sides completes a harmonic pencil centered at it which is transversal to $\Q$, see Figure~\ref{Curva_de_Armonia}.a. These new lines are the \emph{tangent} lines to $\C_\Q$ at the vertices and will be denoted by the corresponding lower case letter. The harmonic bundle of the quadrilateral $a, c, b, d$ is called the \emph{tangent bundle} of $\C_\Q$ and will be denoted $\C_\Q^*$.

Now consider a point $Z\in \C_\Q$ different from the vertices, we call it  \emph{generic}, see Figure~\ref{Curva_de_Armonia}.b. By definition, the four lines from $Z$ to the vertices are a harmonic pencil centered at $Z$. Let $X=(A\vee B)\wedge (C\vee Z)$ and $Y=(A\vee B)\wedge (D\vee Z)$. Then $A, X, B, Y$ is a harmonic set. Therefore, we can recover $Z$ from $X\in A\vee B$ by defining 
\begin{equation}\label{HC-construction}
Y=X\cdot\rho_{A,B} \quad\text{and}\quad Z=(C\vee X)\wedge(D\vee Y)\,.
\end{equation}
But this makes sense for $X$ ranging over all of $A\vee B$ and gives the four vertices (at a harmonic set), so that $\C_\Q$ is parametrized by  $X\in (A\vee B)$ via this construction. We call this the \emph{HC-construction}.

\begin{lemma}[Duality lemma]\label{Lemma_Dualidad}
Points in the harmonic curve $\C_\Q$ are \emph{paired} (i.e., in bijective correspondence) by incidence with lines in its tangent bundle $\C_\Q^*$.
\end{lemma}

\begin{figure}[H]
         \centerline{
          \href{https://arquimedes.matem.unam.mx/harmonic_curves/viewer.html?view=pg3/Duality_between_harmonic_curves_and_bundles.pg3}{\includegraphics[width=7.0cm]{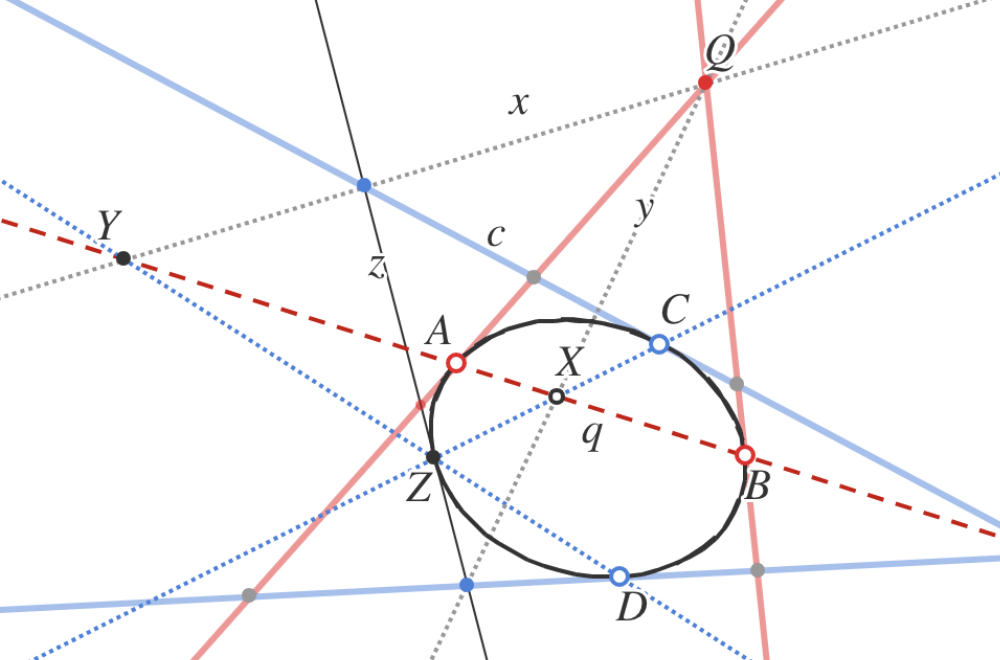}} }  
          \caption{\label{Dualidad_curvas_haces} Incidence of points in a harmonic curve and lines in its tangent bundle.}  
\end{figure}

\proof Let us continue with the notation above, so that $a, c, b, d$, is the quadrilateral whose harmonic bundle is $\C_\Q^*$. As before, these four \emph{generating} lines belong to the bundle because the vertex to which they are tangent (called their \emph{contact point}) can be obtained as the harmonic fourth of their intersection to the other three lines (see Figure~\ref{Curva_de_Armonia}.a). Going further on the {\emph{HC}-construction} \eqref {HC-construction}, and dualizing  it (see Figure~\ref{Dualidad_curvas_haces}): let $Q=a\wedge b$, $x= Q\vee Y$ and $y= Q\vee X$, so that $a, x, b, y$ is generically a harmonic pencil centered at $Q$. Then, $z=(c\wedge x)\vee(d\wedge y)$ is a line of the bundle $\C_\Q^*$, and any such line is uniquely expressed in this way.  

To prove that $Z\in z$ for $X$ different from $A$ and $B$, define $q=A\vee B$ and consider the triangle $QXY$ with respective opposite sides $qxy$. By the definitions, we have $D=C\cdot\rho_{Q,q}$ ($d=c\cdot\rho_{Q,q}$), and by Klein's Triangle Lemma, 
$\rho_{X,x}=\rho_{Q,q}\cdot\rho_{Y,y}$,  then $C\cdot\rho_{X,x}=D\cdot\rho_{Y,y}$ ($c\cdot\rho_{X,x}=d\cdot\rho_{Y,y}$). But $C\cdot\rho_{X,x}\in C\vee X$ and $D\cdot\rho_{Y,y}\in D\vee Y$, so that $C\cdot\rho_{X,x}=(C\vee X)\wedge(D\vee Y)=Z$ (dually, $c\cdot\rho_{X,x}=z$). Hence, the fact that $C\in c$, implies that $Z\in z$ as we wished. $Z$ is called the \emph{contact point} of $z\in\C_\Q^*$ which is the \emph{tangent line} to $\C_\Q$ at $Z$.      \qed

As a corollary, we can express the harmonic curve $\C_\Q$ as a family of harmonic reflections applied to a single point
\begin{equation}\label{A-construction}
\C_\Q=\{\,C\cdot\rho_{X,x}\,|\,X\in (A\vee B)\setminus\{A,B\}\,\}\cup\{A,B\}\,,
\end{equation}
where $x=(X\cdot\rho_{A,B})\vee(a\wedge b)$, which only depends on three points $A, C, B$ and the two tangent lines $a, b$ incident to $A, B$, respectively; we will call this, the \emph{A-construction}.

Clearly, harmonic curves are sent to harmonic curves under projections because projections preserve harmonicity. So that the fact that the classic conic sections are harmonic curves follows from the fact that a circle is a harmonic curve. Indeed, consider an inscribed square to a circle as  a quadrangle. Use the inscribed angle theorem to see that each point in the circle is the center of a transversal pencil to the quadrangle with consecutive lines at angles ${\pi}\over{4}$, which is a harmonic pencil.

\section{Polarities and hyperbolic geometry}

A \emph{polarity} in the plane (in space) is a bijective correspondence between points and lines (planes) that preserves incidence; the terms \emph{polar} of a point, \emph{pole} of a line (plane) or a \emph{polar pair} are used\footnote{An extra hypothesis is required in \cite{CoxeterPG}. Namely, that for some line, the map to the line pencil of its pole be a projectivity. But we do not need to stress this issue.}.
 
\begin{theorem}[Polarity]\label{teo:Polaridad2D}
A harmonic curve $\mathcal{C}$ induces a polarity (expressed by upper and lower case of the same letter) satisfying:
\begin{itemize}
\item[i)]{ $P\in \mathcal{C} \Leftrightarrow  P\in p\,.$}
\item[ii)]{If $P\not\in \mathcal{C}$ then the harmonic reflection $\rho_{P,p}$, with $P$ as center and its non-incident polar line $p$ as mirror, leaves $\mathcal{C}$ invariant.}
\end{itemize}
\end{theorem}

We have already seen a part of item (i) as Lemma~\ref{Lemma_Dualidad} because tangent lines to points in a harmonic curve are defined as their polar lines. The rest of the proof will be given in the next section as a consequence of an analogue in 3D, Theorem~\ref{teo:Polaridad3D}. For the moment, let us make two remarks about this theorem and, assuming it is true, explore some of its profound consequences.

The two mathematicians that father this theorem are Jean-Victor Poncelet (1788-1867) and Karl G. C. von Staudt (1798-1867). Poncelet proved the relation of poles and polars for conic sections using harmonicity (in its metric version), and soon after, von Staudt developed polarities as a general concept and used it as an alternative way to define conic curves within projective geometry with no  metric or algebraic considerations, \cite{VonStaudt}. This definition via polarities is the one Coxeter uses in his influential book \cite{CoxeterPG}, and calls it ``extraordinarily natural and symmetrical'' because it has duality built into it. In general, there are two types of polarities: \emph{euclidian} in which no point is incident with its polar line, and \emph{hyperbolic} when there exist pole and polar incident pairs. The terms used are related to the groups generated by harmonic reflections of non-incident polar pairs. So that von Staudt's definition of a conic curve is equivalent to item (i) of the theorem for a hyperbolic polarity, while Poncelet's results can be rephrased as item (ii). 

As examples of polar pairs, we have named lines and points in Figure~\ref{Dualidad_curvas_haces} according to the upper and lower case rule for poles and polars with respect to the displayed harmonic curve $\C_\Q$. Indeed, a point $Z\in \C_\Q$ and the corresponding line in  $z\in\C_\Q^*$ described in Lemma~\ref{Lemma_Dualidad} constitute a polar pair satisfying (i). 

We now prove that von Staudt's definition of conic curves with mild extra hypothesis gives harmonic curves.

\begin{lemma}\label{vS_gives_HC}
Given a  polarity in the plane, let $\C$ be the set of points that are incident to their polar line and suppose item (ii) of  Theorem~\ref{teo:Polaridad2D} holds. If every line  meets $\C$ in at most two points and $\C$ contains at least three points, then $\C$ is a harmonic curve.
\end{lemma}

\proof Let $A, B, C \in\C$ be three points. By the hypothesis on the lines, they are not collinear.  
Let $a, b$ be the respective polar lines of $A, B$, so that $A\in a$ and $B\in b$. Let $Q=a\meet b$; it is the pole of $q=A\vee B$ because polarities preserve incidence, which also implies that $Q\not\in q$. Finally, let $\Q$ be the quadrangle $A, C, B, D=C\cdot\rho_{Q,q}$. To conclude the proof we show that $\C=\C_\Q$. 

Given $X\in q\setminus \{\,A, B\,\}$, its polar, $x$, is a line through $Q$ different from $a$ and $b$. Let $Y=x\meet q$. Then, since $\rho_{X,x}$ leaves $q$ and $\C$ invariant and $q\cap\C=\{\,A, B\,\}$, it transposes $A$ and $B$, so that $X, A, Y, B$ is a harmonic set. Since the polarity satisfies (ii) of Theorem~\ref{teo:Polaridad2D}, $C\cdot\rho_{X,x}\in\C$, so that the A-construction \eqref{A-construction} for $\C_\Q$ implies that $\C_Q\subset\C$. Finally, given $Z\in\C$ different from $A,B,C$, let $X=(Z\vee C)\meet q$, then $Z=C\cdot\rho_{X,x}$ because the line $Z\vee C$ has no point in $\C$ other than $Z$ and $C$ by hypothesis. Therefore, $\C_Q=\C$.
\qed

As a corollary to the proof of this lemma (and assuming the Polarity Theorem) we can now say which quadrangles $\Q$  inscribed in a harmonic curve $\C$ generate it as its harmonic curve, that is, are such that $\C_\Q=\C$. Precisely \emph{the ones in which the pole of one diagonal line lies in the other diagonal}. And moreover, for any triplet of points $A, C, B$ (with $C$ distinguished) in $\C$, we have such a quadrangle $A, C, B, D=C\cdot\rho_{(a\wedge b), (A\vee B)}$.

The abstract resemblance of these facts to the harmonic fourth construction and the existence of harmonic sets in the projective line lead to the theorem stated and proved bellow. But we must also remark that it is deeply related to the projective model of the hyperbolic plane due to Beltrami and Klein, which he used as an important example for his Erlangen Program, \cite{Erlangen}. 

If we fix a harmonic curve $\C$, any two points $A, B$ in $\C$ define a \emph{hyperbolic line} $q=A\vee B$ and a \emph{hyperbolic reflection} $\eta_q=\rho_{Q,q}$ where $Q$ is the pole of $q$. 

Then we have a subgroup $\Hyp(2)$ of $\Har(2)$, called the \emph{hyperbolic group}, which is generated by all the hyperbolic reflections $\eta_q$ on hyperbolic lines. This group acts on the \emph{inside} of the harmonic curve $\C$ which is taken as the \emph{hyperbolic plane} for the Beltrami-Klein model, and hyperbolic geometry can be built from there, but we have no space here to follow it through; for such treatments see, e.g., \cite{CoxeterNEG, Santalo, GeomVis} . However, we should remark that the generating quadrangles of $\C$ have as diagonals the pairs of perpendicular hyperbolic lines.

The usual way to state the following theorem is in the context of the upper half plane model of the hyperbolic plane and referring to matrix groups. 

\begin{theorem} $\,\, \Hyp(2)\cong\Har(1)$.
\label{thm:HypHar}
\end{theorem}

\begin{figure}[H]
         \centerline{
           \href{https://arquimedes.matem.unam.mx/harmonic_curves/viewer.html?view=pg3/Tangential_map.pg3}
          { \begin{overpic}[abs,width=6.1cm]{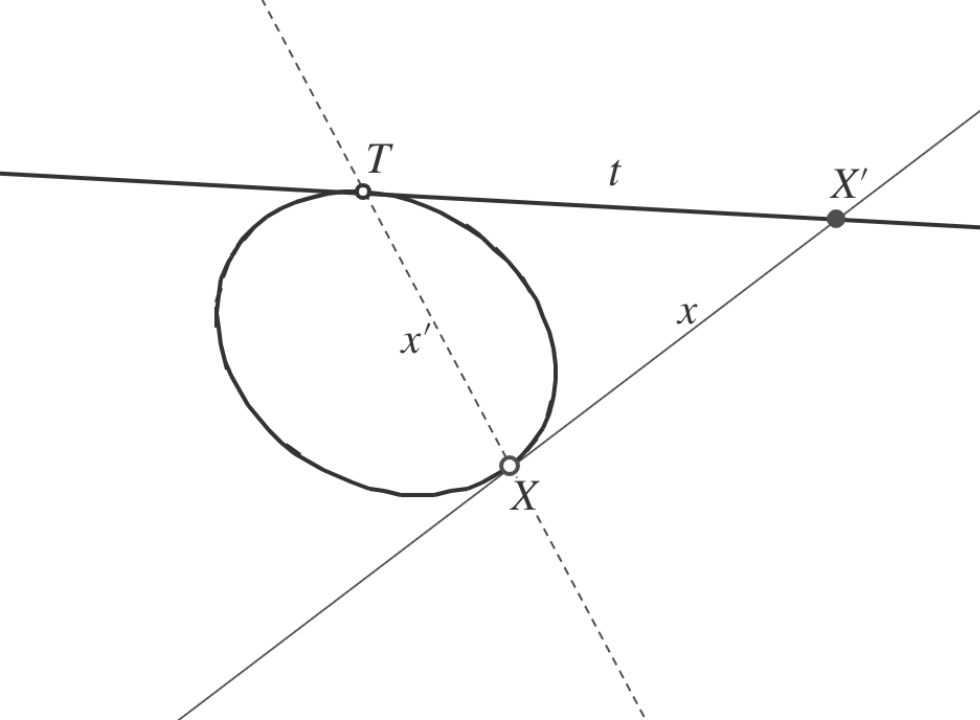}
        \put(-5,0){\Fa}
        \end{overpic} }                 
          \hspace{0.2cm}
     \href{https://arquimedes.matem.unam.mx/harmonic_curves/viewer.html?view=pg3/Generating_quadrangle.pg3} 
         { \begin{overpic}[abs,width=6.1cm]{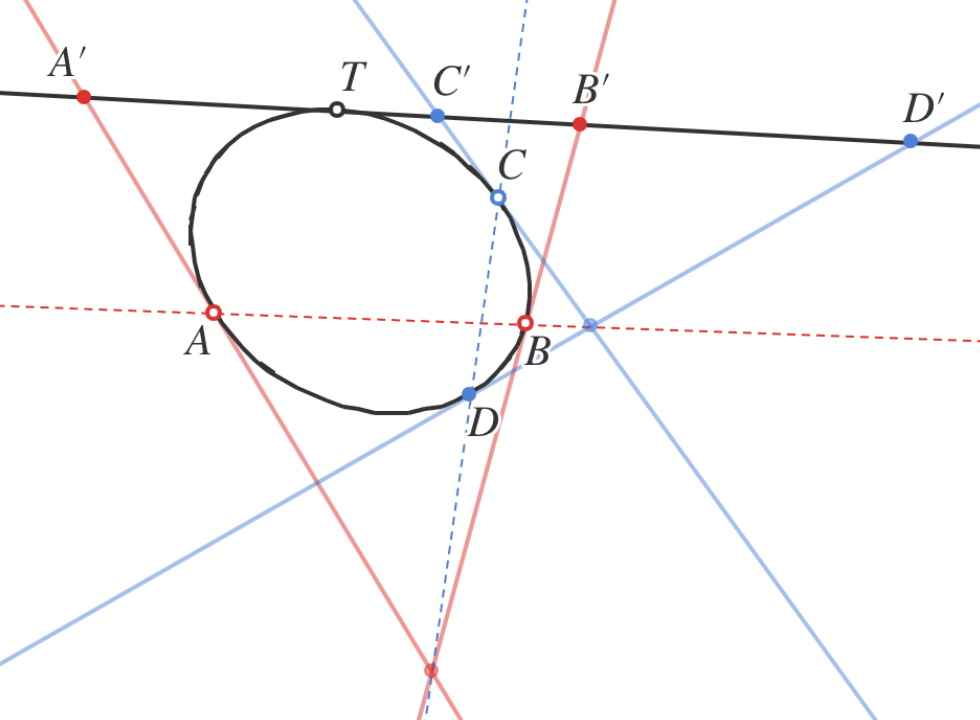}  
              \put(-5,0){\Fb}
           \end{overpic} }}
          \caption{\label{Corr_Tang} {\bf a)} Tangential map. {\bf b)} A generating quadrangle in $\C$ and its corresponding harmonic set in $t$.}  
\end{figure}

\proof
First, we define the \emph{tangential map} from $\C$ to a tangent line, Figure~\ref{Corr_Tang}.a. Let $T$ be a point in $\C$ and let $t$ be its tangent (or polar) line. For every $X\in\C$ other than $T$, let $X^\prime=t\wedge x$, where $x$ is the tangent line to $\C$ at $X$. Taking $T=T^\prime$,  this gives a bijective map $X\leftrightarrow X^\prime$ between $\C$ and $t$, because $x^\prime$ (the polar of $X^\prime\in t$) cuts $\C$ in $T$ and $X$ for $X^\prime\neq T$, see Figure~\ref{Corr_Tang}.a.

Considering $t$ as $\P^1$, the theorem follows from the fact that generating quadrangles of $\C$ and harmonic sets of $t$ correspond to each other, because then harmonic reflections and hyperbolic reflections (the generators of the groups) correspond under the tangential map. 

Let $A, C, B, D$ be a generating quadrangle of $\C$, and let $a, c, b, d$ be their respective tangent lines. By Lemma~\ref{Lemma_Dualidad}, the harmonic bundle of this quadrilateral is the tangent bundle of $\C$ and it contains $t$. Therefore, by the definition of harmonic bundles, we have that $A^\prime, C^\prime, B^\prime, D^\prime$ is a harmonic set, see Figure~\ref{Corr_Tang}.b.   \qed

This proof was the primal motivation for our definition of \emph{harmonic curves}. 

\section{Doubly ruled surfaces}

The following proof of the Polarity Theorem (\ref{teo:Polaridad2D}) is inspired by Dandelin's proof of Pascal's Hexagonal Theorem. Given a conic curve, Dandelin constructs, in \cite{Dandelin}, a hyperboloid of revolution that has it as a plane section; then, using the fact that this surface is doubly-ruled, he obtains a configuration of 6 lines in three dimensional space associated to the six points of the hexagon in the conic, and argues with the geometric-combinatorial properties of the configuration to conclude the proof.  We use the same general idea and get to the same configuration of 6 lines, but instead of hyperboloids of revolution we can now use general ruled surfaces following Hilbert and Cohn-Vossen's construction of ruled surfaces in \cite{HCV}, which appeared in print almost a century after Dandelin's proof, and made clear that they can be constructed by simple incidence arguments.

Consider two lines $a$ and $b$ in three dimensional projective space. They touch if and only if they are coplanar. If this is not the case, they can be called a \emph{generating} pair because for any point $X$ not in them, there exists a unique line through $X$ \emph{transversal} (i.e., with a common point) to $a$ and $b$; namely: 
$$(X\vee a)\wedge(X\vee b)\,.$$
Now consider three lines $a$, $b$, $c$ in \emph{general position} (i.e., each pair is generating, or equivalently, no pair of them is coplanar). The \emph{transversal ruling} to $a, b, c$, denoted $\R(a, b, c)$, 
is the set of lines that are  \emph{transversal} to them (i.e., that touch all three); any such set of lines will be called a \emph{ruling} and its elements are called its \emph{rules}, see Figure~\ref{Reglados}.a. If we denote $\R=\R(a, b, c)$, the above observation implies that $\R$ is parametrized by incidence with the points in any of the three generating lines (through any point in them there passes a unique rule). It will be important to note that, dually, $\R$  is also parametrized by planes containing one of the lines; if we denote planes by greek letters (points and lines  are, respectivelly, upper and lower case latin) we have, for example, that
\begin{equation}\label{eq:ParamPorPlanos}
\mathcal{R}(a,b,c)=\{\,(b\wedge \alpha)\vee(c\wedge\alpha)\,|\,a\subset\alpha\,\}\,.
\end{equation}

\begin{figure}[h]
         \centerline{ 
         \href{https://arquimedes.matem.unam.mx/harmonic_curves/viewer.html?view=pg3/Transversal_ruling_by_blue_lines.pg3}
     { \begin{overpic}[abs,width=6.1cm]{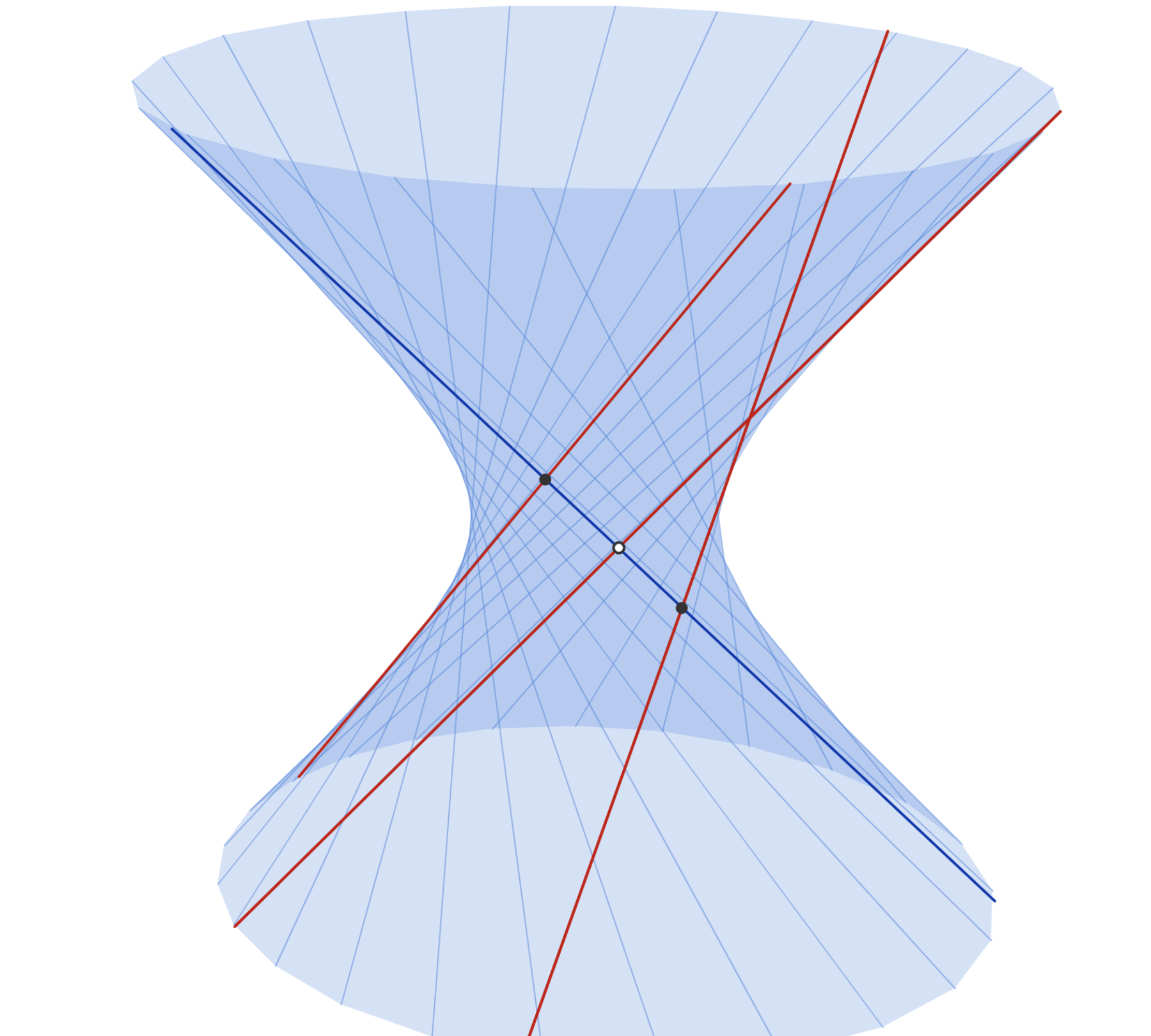}
        \put(-5,0){\Fa}
        \end{overpic} }                 
          \hspace{0.2cm}
         \href{https://arquimedes.matem.unam.mx/harmonic_curves/viewer.html?view=pg3/Transversal_ruling_by_red_lines.pg3} 
         { \begin{overpic}[abs,width=6.1cm]{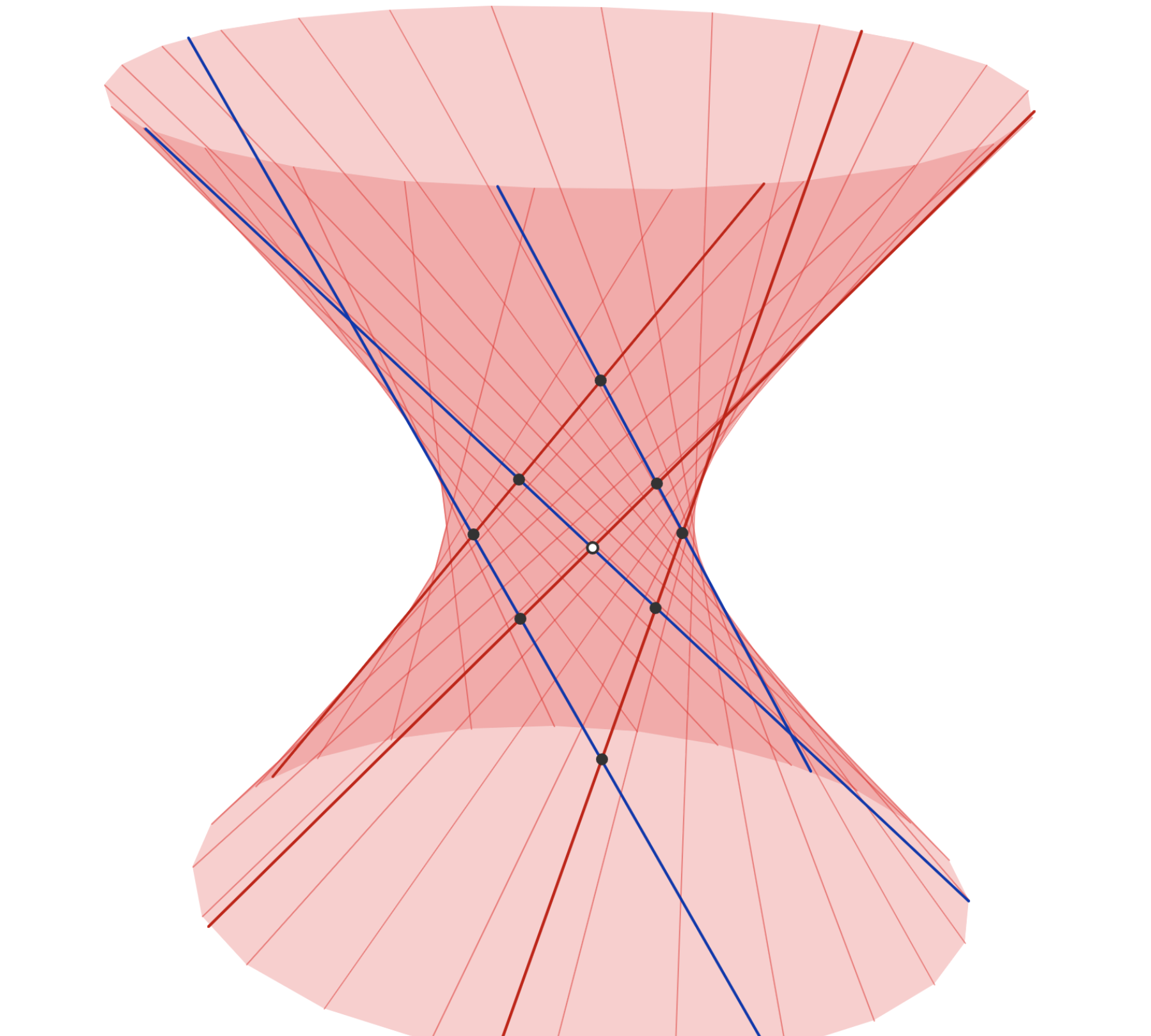}   
           \put(-5,0){\Fb}
         \end{overpic} }}
          \caption{ {\bf a)} The transversal ruling by blue lines to three red lines. {\bf b)} The transversal ruling to any three blue rules contains the three original red lines.} 
\label{Reglados}           
\end{figure}

Every pair of rules in $\mathcal{R}$ is generating, otherwise, their three transversal lines $a, b, c$ would be coplanar. Thus, for any triplet $a^\prime, b^\prime, c^\prime\in\mathcal{R}$ we get a transversal rulling $
\mathcal{R}(a^\prime, b^\prime, c^\prime)$ that contains the original three lines, $a, b, c$; this ruling is an \emph{extension} of $a, b, c$, see Figure~\ref{Reglados}.b. In \emph{real} projective space there is only one extension to a ruling of three lines in general position, but this is not true in all projectives spaces. In what follows we will assume that our projective space does have this property, that is, we will assume it as an axiom:

\medskip
\noindent{\bf Equipal Axiom.}\footnote{Equipal is a classic mexican style of furniture that uses double rulings for bases, \cite{GeomVis}.}
\emph{Three lines in general position belong to a unique ruling.}
\medskip

Later on we will prove it is equivalent to Pappus' Theorem and thus to other classic statements that are commonly adopted as axioms in projective geometry. Another name for it could be the ``Double-ruling Axiom'' because it immediately implies that rulings are matched or paired: any ruling has an \emph{opposite ruling}  
which is the transversal ruling to any three of its rules. The \emph{doubly-ruled surface} (we also refer to it simply as a \emph{ruled surface}) obtained as the union of the rules in any such ruling is also the union of the rules in its opposite ruling.  

Hence, every point on a ruled surface has a \emph{tangent plane}, the one generated by the unique rules through the point in the opposite rulings of the surface.  

\begin{theorem}[Polarity of ruled surfaces]\label{teo:Polaridad3D}
The pairing of points in a ruled surface $\S$ with their tangent planes extends to a polarity of projective space.  Furthermore, if $P\not\in\mathcal{S}$ then $P$ is not incident with its polar plane $\pi$ and the harmonic reflection $\rho_{P,\pi}$, with $P$ as center and $\pi$ as mirror, leaves $\S$ invariant.
\end{theorem}

\proof The ruled surface $\S$ has two opposite rulings $\R$ and $\R^\prime$ such that 
$$\mathcal{S}=\bigcup_{x\in\mathcal{R}}x=\bigcup_{y\in\mathcal{R}^\prime}y\,.$$

To define the polarity induced by $\S$ in its complement, fix three rules $a, b, c$ in the ruling $\R$, and beware that we have inverted the notational use of primes: their transversal ruling is now $\R^\prime=\R(a,b,c)$.

Consider a point $P\not\in\S$; dually, we could start with a non-tangent plane.  

Let $\alpha=a\vee P$.
There is a well defined rule $a^\prime\in \R^\prime$ for which $P\in a\vee a^\prime=\alpha$ (namely, 
$a^\prime=(b\wedge \alpha)\vee(c\wedge\alpha)$, as in \eqref{eq:ParamPorPlanos}). 
Let $A=a\wedge a^\prime\in\S$. Observe that $A$ must be in the polar plane of $P$ because polarities preserve incidence and $P$ is in the polar plane of $A$.

Analougously, we obtain $b^\prime,  c^\prime\in\R^\prime$, for which $P\in b\vee b^\prime=\beta$ and $P\in c\vee c^\prime=\gamma$. Let $B=b\wedge b^\prime$ and $C=c\wedge c^\prime$,  
so that the polar plane to $P$ has to be
$$\pi=A\vee B\vee C\,.$$

\begin{figure}[h]
         \centerline{ 
         \href{https://arquimedes.matem.unam.mx/harmonic_curves/viewer.html?view=pg3/Polar_plane_of_a_point_with_respect_to_a_ruled_surface.pg3}
     { \begin{overpic}[abs,width=6.1cm]{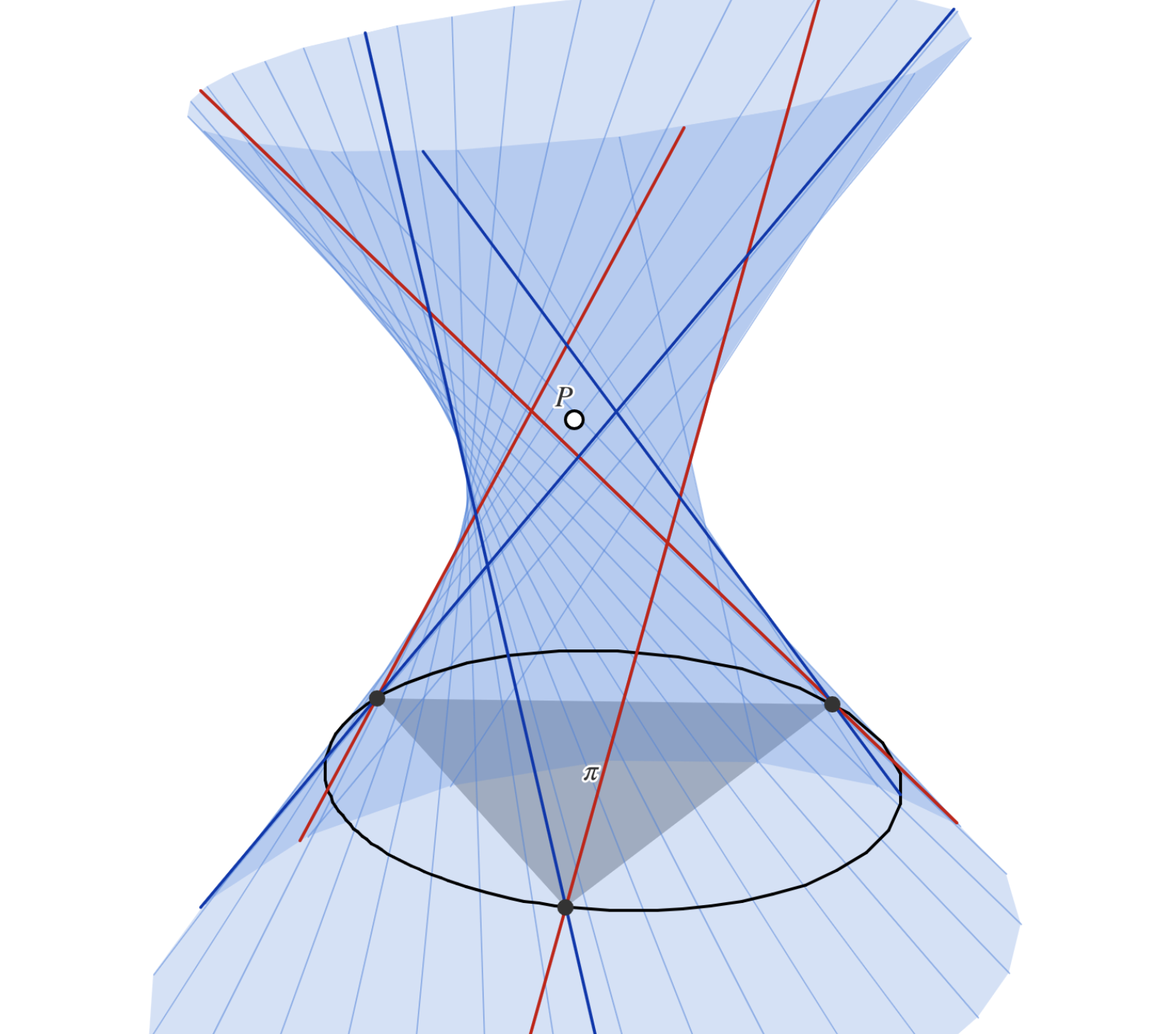}
        \put(-5,0){\Fa}
        \end{overpic} }                 
          \hspace{0.2cm}
         \href{https://arquimedes.matem.unam.mx/harmonic_curves/viewer.html?view=pg3/Dandelin_configuration_from_ruled_surface.pg3} 
         { \begin{overpic}[abs,width=6.1cm]{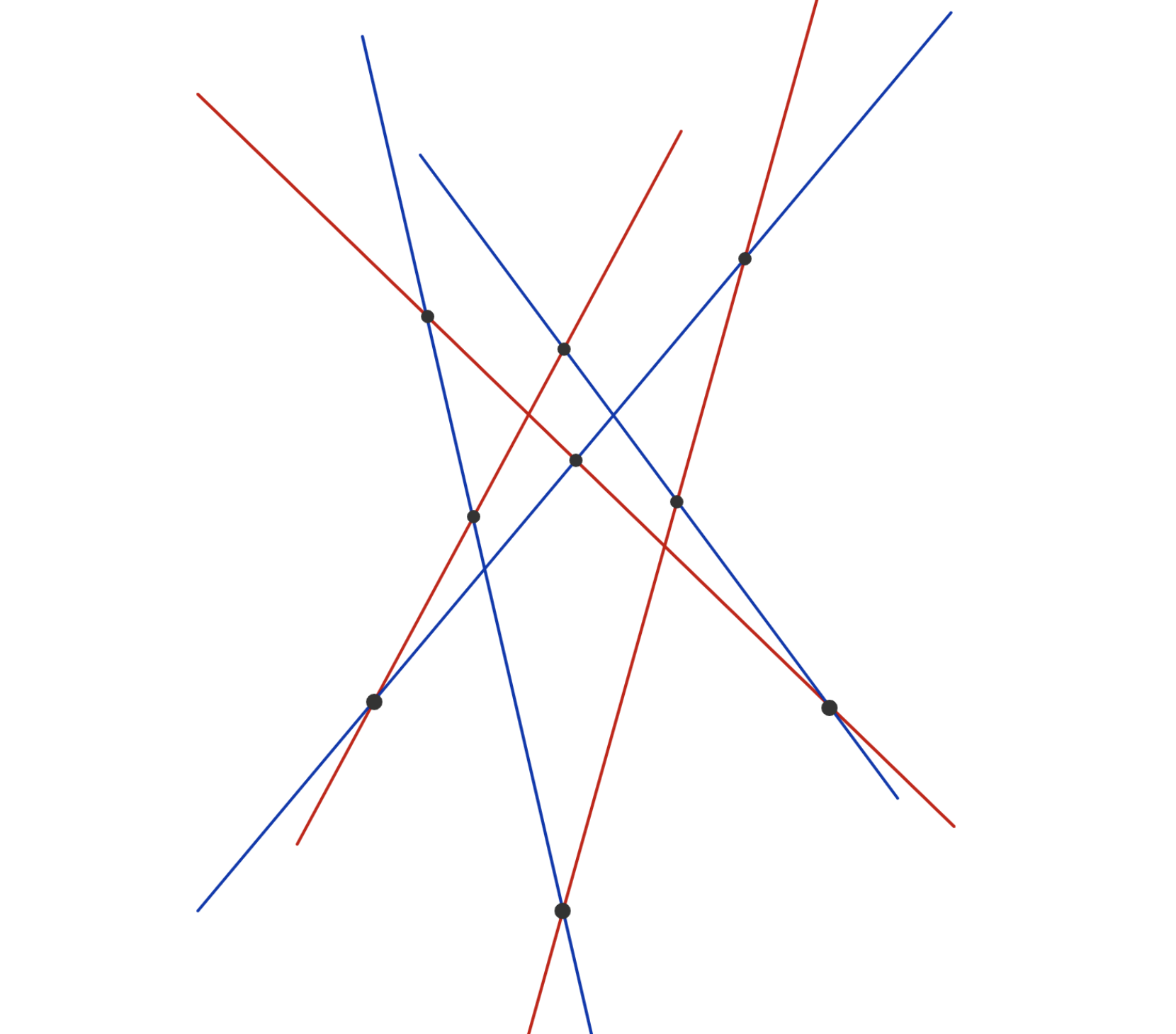}   
           \put(-5,0){\Fb}
         \end{overpic} }}
          \caption{  {\bf a)} Polar plane $\pi$ of a point $P\not\in\S$ with respect to a ruled surface $\S$. {\bf b)} Dandelin configuration.} 
\label{Reglados}           
\end{figure}

If we had started, dually, with a non tangent plane $\pi$ we would have found $P$ as the intersection of the three tangent planes at $A=a\wedge\pi$, $B=b\wedge\pi$,  $C=c\wedge\pi$; and $a^\prime, b^\prime, c^\prime$ would be the rules in $\R^\prime$ passing through $A, B, C$ respectively. So that the pairing of points and planes is now well defined.

We have distinguished what we will call a \emph{Dandelin configuration}: six lines of two \emph{types} or \emph{colors}, three of each, $a, b, c$ and $a^\prime, b^\prime, c^\prime$---unprimed and primed in the text,  
red and blue in the pictures as in Figure~\ref{Reglados}.b--- such that a pair of them touch if and only if they have opposite types. This produces nine \emph{basic} points and nine \emph{tangent} planes by the ``meet'' ($\wedge$) or ``join'' ($\vee$) of lines of different colors; but it also comes with a derived configuration of other lines and planes that naturally arise from them. The geometric richness of this configuration, closely related to the combinatorics of $3\times 3$ determinants, is what Dandelin exploited in \cite{Dandelin}; and we follow suit.

Now, we will prove that the harmonic reflection, $\rho_{P,\pi}$, with center $P$ and mirror $\pi$ interchanges the lines $a, b, c$ respectively with $a^\prime, b^\prime, c^\prime$ in the opposite ruling. By the triangular symmetry of the construction, it will sufffice to prove that:

\begin{quote} $\bullet$ \emph{in the tangent plane to $A$, ${\alpha=a\vee a^\prime}$, the lines $a, A\vee P, a^\prime,\alpha\wedge \pi$  are a harmonic pencil centered at $A$.}\end{quote} 

Because this happens if and only if $\rho_{P,\pi}$ interchanges the lines $a$ and $a^\prime$. 

\begin{figure}[b]
          \centerline{
\href{https://arquimedes.matem.unam.mx/harmonic_curves/viewer.html?view=pg3/A_Dandelin_configuration_from_a_harmonic_curve_with_pencil.pg3}{\includegraphics[width=9cm]{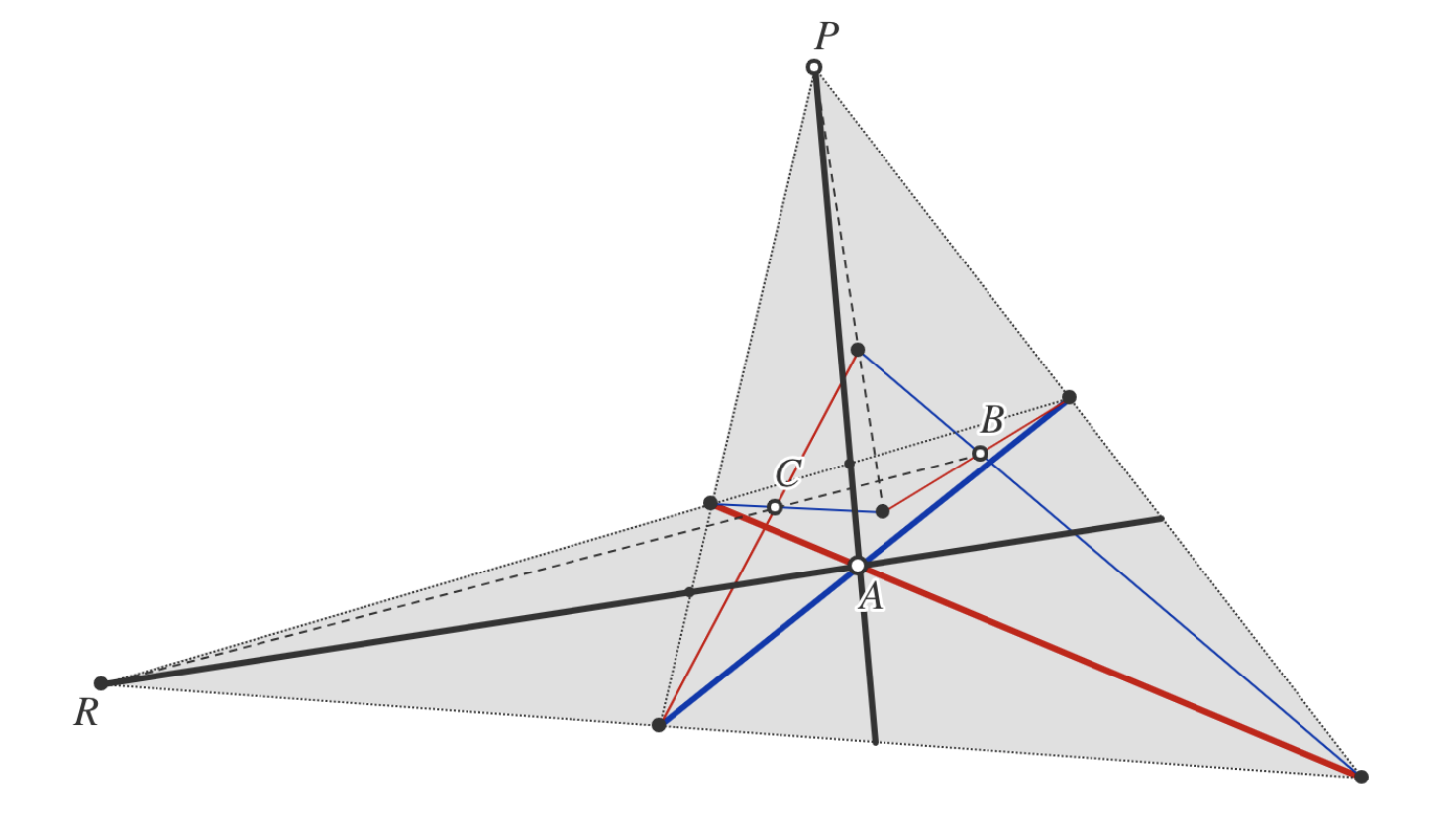}}}
                   \caption{\small The Dandelin configuration given by the point $P\not\in\S$                    with the harmonic pencil $a, A\vee P, a^\prime, A\vee R= \alpha\wedge\pi$ in the plane $\alpha=a\vee a^\prime$.} 
\label{fig:Conf_Dandelin} 
\end{figure}

The tangent plane $\alpha=a\vee a^\prime$ contains five of the nine basic points of our Dandelin configuration. Namely, the $\alpha$-quadrangle:
$$a\wedge b^\prime, b\wedge a^\prime, a\wedge c^\prime, c\wedge a^\prime\,,$$ 
with its center $A$ and its diagonals are $a$ and $a^\prime$. The remaining four basic points outside of $\alpha$, break naturally into two pairs whose generated lines are incident with the two diagonal points of the $\alpha$-quadrangle. This follows because these diagonal points can be seen as the intersection of three tangent planes. Namely, ${P=\alpha\wedge\beta\wedge\gamma}$ and $R=\alpha\wedge(b\vee c^\prime)\wedge(c\vee b^\prime)=\alpha\wedge(B\vee C)\in\alpha\wedge\pi$ (see Figure~\ref{fig:Conf_Dandelin}).

Thus, $\rho_{P,\pi}$ interchanges the rules $a$ and $a^\prime$. Analogously, it interchanges $b$ with $b^\prime$ and $c$ with $c^\prime$. Then, it gives a bijection between the transversal rulings of $a, b, c$ and $a^\prime, b^\prime, c^\prime$, which are $\R^\prime$ and $\R$ respectively,  because a line transversal to  $a, b, c$ is sent by $\rho_{P,\pi}$ to a line transversal to $a^\prime, b^\prime, c^\prime$ and viceversa. Therefore, $\rho_{P,\pi}$ leaves $\mathcal{S}$ invariant, as we wished to prove.

In particular, since a harmonic reflection sends a line to a line concurrent with the mirror and coplanar with the center, 
our definition of the polarity does not depend on the choice of generating rules $a, b, c$.

Finally, the proof that the polarity we have defined preserves incidence follows in cases, but in a straightforward manner from the fact that if the tangent plane to a point in $\S$, say $A$ as above, contains a point not in $\S$, say $P$, then the polar plane of $P$ contains $A$.
\qed

Observe that, because of the incidence invariance, the polarity extends naturally to a pairing of lines. The polar of a line $\ell$ is the intersection of all the polar planes of its points, or of any two of them.

This polarity theorem asserts that what one sees as the contour of a ruled surface is exactly its section with the polar plane of the viewpoint. 
Sections and the contour of projections coincide.  
We now prove that sections of ruled surfaces are harmonic curves, and that the corresponding harmonic bundle is the projection from the pole of any one of the two rulings.  

\smallskip
\noindent\emph{Proof of Theorem~\ref{teo:Polaridad2D}.}  Consider a harmonic curve, $\mathcal{C}$, in a plane $\pi$. Our basic aim is to prove that 

\begin{quote} $\bullet$ \emph{there exists a ruled surface $\S$ that has $\C$ as a section,}
\end{quote}

that is, such that $\C=\S\cap\pi$. This will induce the desired polarity in $\pi$ to complete the proof of the theorem. 

By definition, $\C$ is the harmonic curve of a quadrangle $A, C, B, D$. Let $a$ and $b$ be the tangents at $A$ and $B$, respectively; and let $Q=a\wedge b$, $q=A\vee B$. We know that $D=C\cdot\rho_{Q,q}$ and that $\C$ is obtained by the \emph{A}-construction \eqref{A-construction}. 

\begin{figure}[h]
          \centerline{
\href{https://arquimedes.matem.unam.mx/harmonic_curves/viewer.html?view=pg3/A_Dandelin_configuration_from_A_construction_of_harmonic_curve.pg3}
     { \begin{overpic}[abs,width=6.1cm]{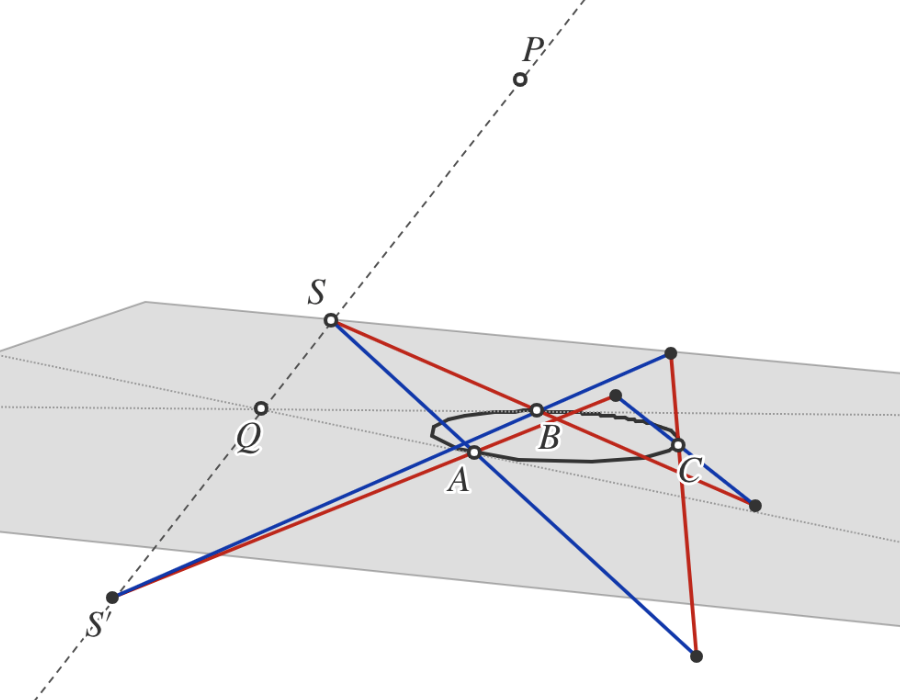}
        \put(-5,0){\Fa}
        \end{overpic} }                 
          \hspace{0.2cm}
       \href{https://arquimedes.matem.unam.mx/harmonic_curves/viewer.html?view=pg3/Ruled_surface_from_A_constr_of_harm_curve.pg3} 
         { \begin{overpic}[abs,width=6.1cm]{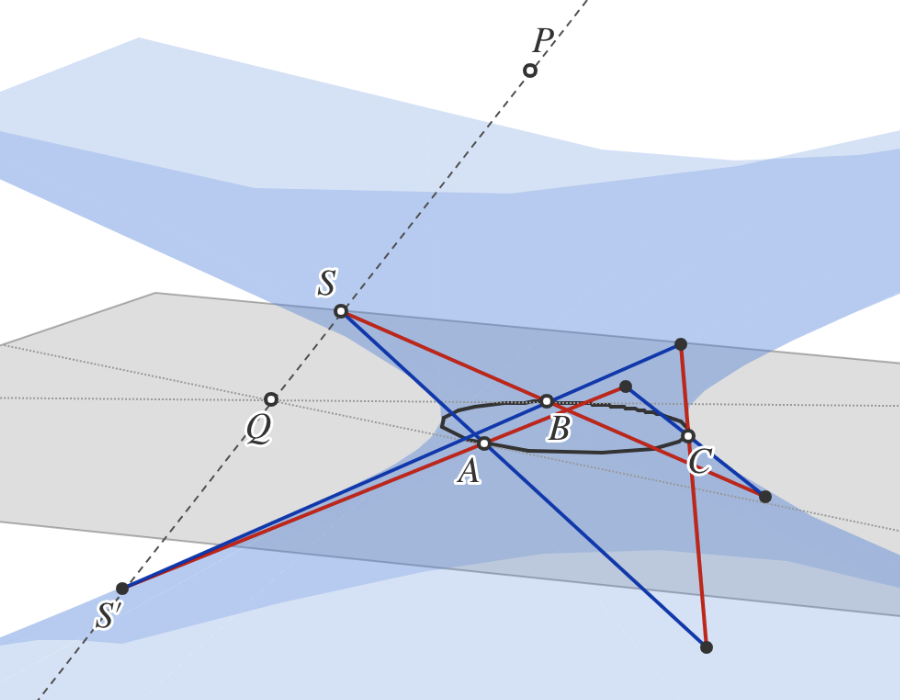}   
           \put(-5,0){\Fb}
         \end{overpic} }}
  \caption{  {\bf a)} A Dandelin configuration arising from the input of the \emph{A}-construction of $\C$ in a plane~$\pi$. {\bf b)} The corresponding ruled surface that intersects $\pi$ in $\C$.\vspace{-6pt}}  
  \label{CA_a_SupReg}
\end{figure}

Choose two points $P$ and $S$ not in $\pi$ and colinear with $Q$ (see Figure~\ref{CA_a_SupReg}.a). 

Let $S^\prime=S\cdot\rho_{P,Q}$. Since $S\neq S^\prime$, the four lines from $S$ and $S^\prime$ to $A$ and $B$ can be colored red and blue so that only lines of opposite colors touch. Finally, consider the red (blue) line through $C$ transversal to the two blue (red) lines. We now have a Dandelin configuration of six lines colored red and blue: let $\S$ be the doubly ruled surface it defines (Figure~\ref{CA_a_SupReg}.b). By construction, $P$ and $\pi$ are a polar pair with respect to $\S$.

The polarity induced by $\mathcal{S}$ restricts naturally to a polarity in the plane $\pi$ as follows. The polar line of a point in $\pi$ is the intersection of its polar plane with $\pi$, and the pole of a line in $\pi$ is the pole of the plane it generates with $P$ ---or the intersection with $\pi$ of its polar line. In particular, item (i) of Theorem~\ref{teo:Polaridad2D} follows for $\mathcal{S}\cap\pi$.

Since harmonic reflections preserve the planes through their center, those for non-incident polar pairs with pole in $\pi$, restrict to harmonic reflections of $\pi$ that leave  $\mathcal{S}\cap\pi$ invariant. 
Therefore, item (ii) of Theorem~\ref{teo:Polaridad2D} follows for $\mathcal{S}\cap\pi$.

That $\C=\mathcal{S}\cap\pi$ now follows from Lemma~\ref{vS_gives_HC} and its proof (identifying $A, B, C$ in both settings) because a line that intersects $\S$ in three different points is easily seen to be a rule of $\S$ and $\pi$ contains no such rules.
 \qed
 
 Observe that, within the above framework, for any point in $\mathcal{S}\cap\pi$ the intersection with $\pi$ of its tangent plane to $\mathcal{S}$ is the projection to $\pi$ from $P$ of any of its two rules. So that we may state the following theorem as a corollary to the preceding proofs.
 
 \begin{theorem}\label{teo:CH-SR}
Harmonic curves are the sections of ruled surfaces with non-tangent planes. Moreover, harmonic bundles are the projection of rulings from external points, and the tangent bundle of a section of a ruled surface is the projection from the corresponding pole of any of its two rulings. \qed
 \end{theorem}
 
Finally, we prove the following theorem, making the appropriate remarks to acknowledge Dandelin's original proof of Pascal's Hexagon Theorem that inspired our treatment.

\begin{theorem}\label{Papus_Equipal}
The Equipal Axiom is equivalent to Pappus' Theorem.
\end{theorem}

\proof First, we must state Pappus' Theorem: 
\begin{quote} 
$\bullet$ \emph{The opposite sides of a planar hexagon whose vertices lie alternatively in two lines, meet in three collinear points.}
\end{quote}

Let $a_0$ and $b_0$ be  coplanar lines with points $B_1, B_2, B_3 \in a_0$ and $A_1, A_2, A_3 \in b_0$, so that the hexagon of Pappus' hypothesis is $A_1, B_2, A_3, B_1, A_2, B_3$ considered  cyclically, and the theorem asserts that the three ``\emph{Pappus' points}''
$$P_i=(A_j\vee B_k)\meet(A_k\vee B_j)\,,$$
where $\{i, j, k\}=\{1,2,3\}$, are collinear. 

The hypothesis of Pascal's Theorem is that the six points named above lie not on two lines, but on a harmonic curve and the conclusion is exactly the same. Dandelin's proof  
considers rules ($a_i$ and $b_i$, $i=1,2,3$) through the vertices alternatively in the two rulings of a ruled surface. For the case of Pascal, this would now follow immediately from Theorem~\ref{teo:CH-SR}; for Pappus, we need to work a little more because the plane $\pi=a_0\vee b_0$ will turn out to be a tangent one.

Let $a_1, a_2$ be a pair of generating lines that meet $\pi$ in $A_1, A_2$ respectively. Let $\R^\prime=\R(a_0, a_1, a_2)$ so that $b_0\in\R^\prime$ and let $b_1, b_2, b_3\in\R^\prime$ be the rules through $B_1, B_2, B_3$ respectively. Now, let $\R=\R(b_0, b_1, b_2)$ so that $a_0, a_1, a_2\in\R$ and finally, let $a_3\in\R$ be the rule through $A_3\in b_0$.

We have defined eight lines of two types or colors, $a_i$ and $b_j$ with $0\leq i, j\leq 3$, such that all pairs of different color except one do meet, namely, $a_i$ meets $b_j$ for all $i\neq3\neq j$. The Equipal Axiom implies that $\R=\R(b_0,b_1,b_3)$ and thus, that $a_3\in\R$ meets $b_3$. But moreover, the Equipal Axiom follows if this is always true for this general setting of eight lines, because it implies $\R=\R(b_0, b_1, b_3)$ letting $a_3$ run in all of $\R$; and then moving the $b_j$'s around $\R^\prime$, this implies that $a_0, a_1, a_2$ extend to the unique ruling $\R$.

\begin{figure}[h]
          \centerline{
          \href{https://arquimedes.matem.unam.mx/harmonic_curves/viewer.html?view=pg3/Pappus_and_Equipal_axioms_equivalence.pg3}{\includegraphics[width=11.0cm]{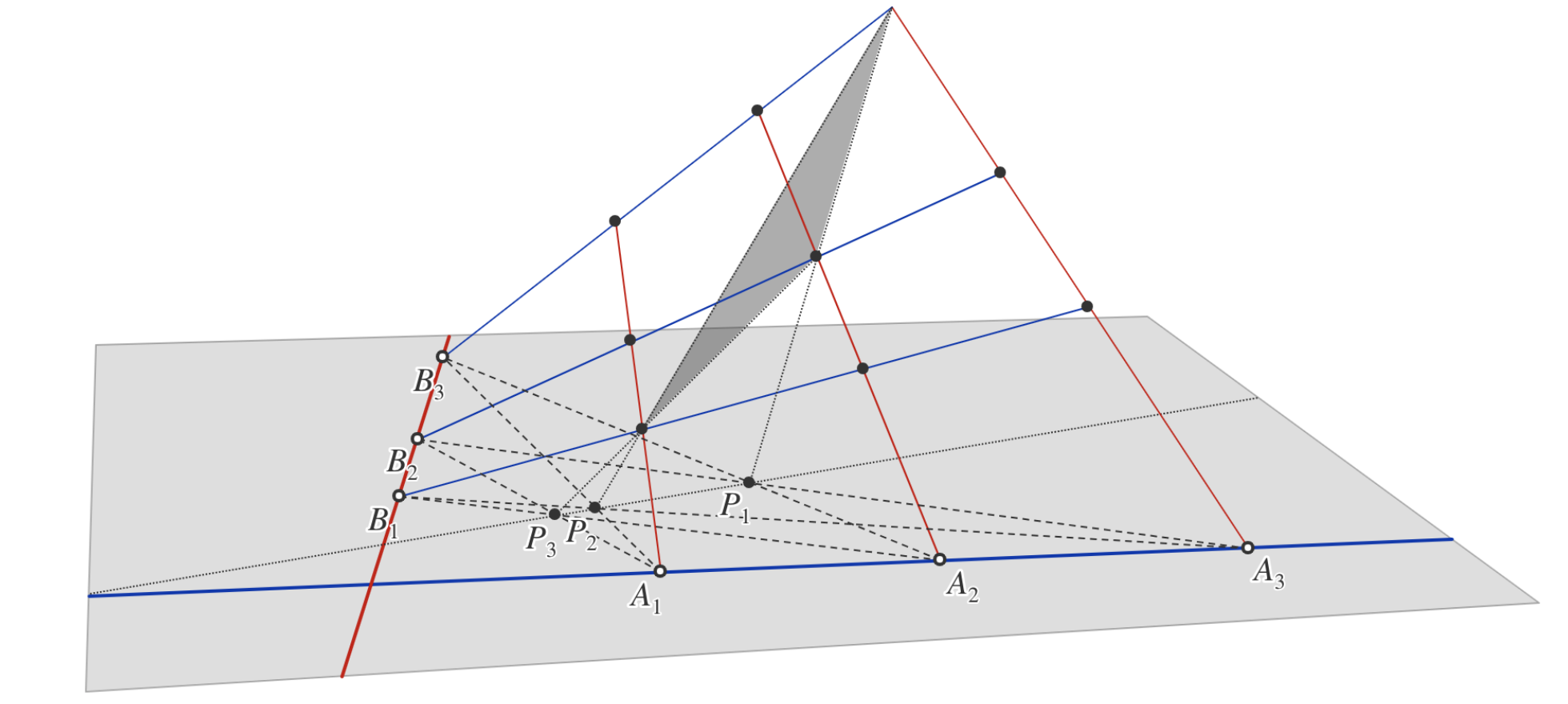} } }  
   \vspace{-6pt}
                   \caption{\small {A Dandelin configuration over a plane with a Pappus configuration.} } 
\label{Pappus-Equipal}
\end{figure}

So, we are left to prove that the Pappus's points $P_1, P_2, P_3$ are collinear if and only if $a_3$ and $b_3$ meet, see Figure~\ref{Pappus-Equipal}.

Suppose that $a_3$ and $b_3$ meet. Then $a_i$ and $b_j$, with $i,j\in\{1,2,3\}$ form a Dandelin configuration. For any such $i,j$ we have that
$$A_i\vee B_j=(a_i\vee b_j)\meet \pi\,.$$

So that the Pappus' points may be seen as lines intersecting $\pi$:
\begin{equation}\label{Pappus-lines}
P_i=((a_j\vee b_k)\meet(a_k\vee b_j))\meet\pi=((a_j\meet b_j)\vee(a_k\meet b_k))\meet\pi\,,
\end{equation} 
for $\{i,j,k\}=\{1,2,3\}$. But these three lines meet pairwise, therefore they lie in a plane that defines the \emph{Pappus' line}:
$$p=((a_1\meet b_1)\vee(a_2\meet b_2)\vee(a_3\meet b_3))\meet\pi\,,$$
which proves Pappus' Theorem and, for non-tangent planes $\pi$ constitutes Dandelin's proof of Pascal's Theorem.

We are left to prove that Pappus implies the Equipal Axiom which, as we have seen, follows from proving that $a_3$ meets $b_3$ assuming that $P_1, P_2, P_3$ lie in a line $p\subset\pi$. Observe that \eqref{Pappus-lines} still holds for $i=3$ (and $\{j,k\}=\{1,2\}$), so that 
$$\delta=p\vee((a_1\meet b_1)\vee(a_2\meet b_2))\,$$
is a plane because the two lines meet at $P_3$.  It contains the lines 
$$\ell_1=P_1\vee(a_2\meet b_2)\quad \text{and}\quad\ell_2=P_2\vee(a_1\meet b_1)\,,$$ 
which give us a point 
$W=\ell_1\meet \ell_2$. To see that $W\in a_3$ and $W\in b_3$ to conclude the proof, observe that $W$ can be seen as the intersection of three planes in two ways; namely,
of $\delta, (a_3\vee b_1), (a_3\vee b_2)$ and $\delta, (a_1\vee b_3), (a_2\vee b_3)$. 
\qed

\section{Loose ends on axioms and projectivities.}

What we have done in this paper works verbatim in the abstract setting of projective geometry. The included images, which are designed to help the reader to understand intuitively the abstract arguments are, of course, drawn in euclidean 2 and 3 dimensional spaces. However no euclidean arguments are used except for a few statements, e.g., when showing that harmonic curves are conic sections and for the model of the hyperbolic plane (the plane not the group!). 
Since Projective Geometry has always been ground for considerations about math foundations, we think it is appropriate to close with a few remarks concerning axiomatics and the rigorous mathematical content of our proofs.

The axioms on which our presentation is based and all its theorems are proved are the following. A \emph{projective space} consists of a ground set, or \emph{space}, of \emph{points} with a well defined family of subsets called \emph{lines}, satisfying:

\begin{enumerate}
\item Any two distinct points $A$ and $B$ lie on a unique line $A \vee B$.
\item If $A, B, C, D$ are four distinct points and lines $A \vee B$ and $C \vee D$ meet, then the lines $A \vee C$ and $B \vee D$ also meet.
\item There are two lines that do not meet. 
\item Lines have more than two points. 
\item The harmonic fourth of three collinear points is neither of them.
\item The Equipal Axiom.
\end{enumerate} 

These axioms are a variation of those commonly used (e.g., the ones suggested by Stillwell in \cite{JS}). The main difference being the replacement of Pappus's Theorem by the Equipal (or Double-rulings) Axiom. Axioms~1 and 2 are the fundamental \emph{incidence} axioms. The statement of Axiom~2 is attributed to Pash and Veblen; it cleverly says that two lines meet if and only if they are \emph{coplanar} without  the need of having planes previously defined. Axiom~3 means the space is at least three dimensional (more than a plane), and it is known to be equivalent to Desargues' Theorem. Axiom~4 is required for geometry to become interesting and not simply some trivialized set theory. 

The two final axioms depend on some further development of the theory; they are not primitive. Axiom~5 guarantees that the \emph{ground field} is of characteristic different from $2$ or, equivalently, that the geometry does not contain the Fano Plane. The \emph{characteristic} of a projective space can be defined geometrically using the harmonic fourth construction; essentially from how far can one go in a \emph{harmonic sequence} without returning. It is needed here to make sense of harmonic curves (and that harmonic reflections are not the identity) because it implies that harmonic sets that have exactly four points do exist. Axiom~6 is a required additional principle for geometry to be rich enough to have a deep relation to other classic branches of mathematics like analysis and topology; bellow, we will discuss the several versions it may adopt. 

From the first 4 axioms, \emph{flats} can be defined as the \emph{closed} subsets under the operation of taking lines, and then, the \emph{dimension} of a flat is obtained as one less than the number of points needed to \emph{generate} it; so that
\emph{planes} are defined as flats of dimension 2, \cite{VeblenYoung}. The incidence properties of planes and lines in a space of dimension 3 are obtained from this; and the Hilbert-Cohn Vossen construction of ruled surfaces follows, making sense of the statement of the Equipal Axiom. 

Since the Equipal Axiom 
is equivalent to Pappus' Theorem, the arithmetization of projective space yields as \emph{ground field} a commutative one. In \cite{JS} commutativity is proved from Pappus' Theorem and the ground field is described from scratch. The necessity of the axiom is proved by constructing a projective space over a non-commutative field like the quaternions.

Another widely used version of Axiom~6 is as the uniqueness part of the Fundamental Theorem of Projective Geometry. It is usually stated in the context of planar geometry where Axioms 1 and 2 become appealingly \emph{dual} (and Axiom 3 is false but replaced by Desargues' Theorem).
A \emph{projectivity} is defined as the composition of \emph{projections} between (points in) lines or (lines in) concurrent pencils; they are always bijections. It is not hard to construct a projectivity determined by it's (arbitrary) effect on (any) three elements of it's domain. This is the existence part of the Fundamental Theorem. However, the uniqueness is proved to be equivalent to Pappus' Theorem, so one must be assumed to prove the other, see \cite{CoxeterPG}, \cite{VeblenYoung}. 

We think that Axiom 3 is natural because, among other things, it responds and gives credit to the original motivation for the creation of Projective Geometry which is renaissance perspective, for which dimension 3 is quintessential. But then, if one thinks about projectivities between non coplanar lines in a three dimensional projective space, one is naturally lead to consider ruled surfaces. Indeed, given a projectivity from a line $a$ to a non coplanar line $b$, the set of lines joining a point in $a$ to its image in $b$ turns out to be a ruling. So, the Equipal Axiom is intimately related to the uniqueness of projectivities given by three arbitrary values (the projectivity is determined by the extended ruling of three lines). Moreover, this association of a set of lines to a map between lines is also a classic idea. It is the dual of how Jakob Steiner (1796–1863) defined conic curves in a purely projective manner; and is a natural, visual way of presenting them, see e.g., \cite{VeblenYoung}. 

Projective geometry is remarkable in many ways. One of them is the importance of some 
mathematical notions that were worked and experimented within it long before their abstract general acknowledgement. For example, projectivities were masterly used almost one century before the notion of sets and the language of abstract functions was stablished; moreover, together they constitute what we now call a groupoid (defined in the mid XX century within category theory). And of course, there is the leading role it played in broadening our notion of geometry and its influence on the dawn of topology. There are many ways to approach it and present it. We hope this paper contributes to the awareness of its cultural significance and the convenience and possibility of finding ``its way down into secondary schools'', \cite{Lehmer}, or at least into early undergraduate courses.

\bmhead{Acknowledgments} J. Bracho acknowledges support from PAPIIT-UNAM Project IN109023.

\section*{Declarations} The authors declare that there is no conflict of interest.

\end{document}